# Stable source reconstruction from a finite number of measurements in the Multi-frequency Inverse Source Problem

Mirza Karamehmedović[1], Adrian Kirkeby[1,2,3], and Kim Knudsen[1]

[1]Department of Applied Mathematics and Computer Science, Technical University of Denmark, Kgs. Lyngby, Denmark
[2]Department of Mathematical Sciences, Norwegian University of Science and Technology, Trondheim, Norway
[3]adrki@dtu.dk

March 13, 2018

## Abstract

We consider the multi-frequency inverse source problem for the scalar Helmholtz equation in the plane. The goal is to reconstruct the source term in the equation from measurements of the solution on a surface outside the support of the source. We study the problem in a certain finite dimensional setting: From measurements made at a finite set of frequencies we uniquely determine and reconstruct sources in a subspace spanned by finitely many Fourier-Bessel functions. Further, we obtain a constructive criterion for identifying a minimal set of measurement frequencies sufficient for reconstruction, and under an additional, mild assumption, the reconstruction method is shown to be stable. Our analysis is based on a singular value decomposition of the source-to-measurement forward operators and the distribution of positive zeros of the Bessel functions of the first kind. The reconstruction method is implemented numerically and our theoretical findings are supported by numerical experiments.

## 1 Introduction

This paper concerns the multi-frequency inverse source problem (ISP) for the scalar Helmholtz equation in the plane. For a compactly supported source, the frequency-dependent radiated field is measured on a curve outside the source support, at several frequencies. The inverse problem is then to estimate the source from these measurements.

Let $D, D_0 \subset \mathbb{R}^2$ be concentric, open discs of radii $R$ and $R_0$ respectively, with $0 < R_0 < R$. Let $s \in L^2(D_0)$ be the source, such that $\mathrm{supp}(s) \subset D_0$. For a fixed, positive angular frequency $\omega$, let $k = \omega/c$ be the corresponding wave number. Here, $c > 0$ is the constant wave speed and is assumed to be known. Throughout the rest of paper, we will refer to $k$ as *frequency*, since $k$ is directly proportional to $\omega$. Let $u_k$ be the unique solution of

$$\begin{cases} (\Delta + k^2)u_k(x) = s(x), & x \in \mathbb{R}^2, \\ \lim_{\|x\|\to\infty} (\partial_{\|x\|} - ik)u_k = 0 & \text{uniformly for } x/\|x\| \in S^1, \end{cases} \qquad (1)$$

given by

$$u_k(x) = \int_{D_0} G_k(x,y) s(y) dy \quad x \in \partial D. \qquad (2)$$

Here, $G_k(x,y) = -\frac{i}{4} H_0^{(1)}(k|x-y|)$ is the Hankel function of the first kind and of order zero, and $S^1$ is the unit circle. At a fixed frequency $k$ we define the forward operator

$$\begin{aligned} F_k \colon L^2(D_0) &\to L^2(\partial D), \\ s &\mapsto U_k = u_k\big|_{\partial D}. \end{aligned} \qquad (3)$$

$F_k$ is a linear, compact operator (see for instance [6] for a proof of compactness). We call $U_k$ the *measurement at frequency $k$*.

Given a set of frequencies $Q = \{k_j\}_{j \in I}$, where $I$ is some index set, and a corresponding set of measurements $\{U_{k_j}\}$, the multi-frequency inverse source problem is to estimate $s \in L^2(D_0)$ such that

$$U_{k_j} = F_{k_j} s \quad \text{for all } k_j \in Q.$$

This problem has been investigated by several authors in recent years, for example in [4, 5, 6, 16, 26]. It has a wide range of applications, for example in optical tomography [3], antenna control [18] and in the acoustic reconstruction part of several hybrid imaging methods, most notably photo- and thermoacoustic tomography [2, 19], see also [1, 6].

The single-frequency ISP, i.e., when one only has a measurement at a single frequency, is not uniquely solvable due to $F_k$ having a large kernel. In fact $\ker(F_k)$ consists of all so-called non-radiating sources (at frequency $k$), including functions like $s = (\Delta + k^2)u; u \in C_c^\infty(D_0)$, see also [9]. Stability estimates for the single-frequency ISP in terms of a description of the singular value spectrum are shown in [16]. For the multi-frequency ISP, uniqueness and stability of the solution are known when measurements are available over a set of frequencies with an accumulation point [6], and for a set of frequencies coinciding with eigenvalues of the Dirichlet-Laplacian [11](both situations necessarily involve infinitely many measurements). In [1] it was shown that, for arbitrary $s \in L^2$, measurements at a finite number of frequencies are insufficient to obtain a unique solution to the ISP. Both iterative and direct numerical inversion methods have been proposed in the literature [7, 11, 12, 26].

In contrast to the above-mentioned references, in this paper we consider a more practical ISP of finite-dimensional nature. We consider measurements at a finite number of frequencies and reconstruction of sources in a certain finite-dimensional subspace spanned by low-index Fourier-Bessel functions. The main result shows that for a few, carefully chosen frequencies, we can uniquely recover such functions, and that reconstructions are robust with respect to noise in the measurements.

The rest of the paper is organized as follows: In Section 2, we introduce the Fourier-Bessel basis and singular value expansion of $F_k$, and show how this motivates the choice of a finite-dimensional source space. In Section 3, we construct a matrix equation relating the multi-frequency measurements and the finite-dimensional source, and show how the invertibility of this equation depends on the measurement frequencies. Building on these results, Section 4 shows how distribution of the zeros of Bessel functions allows us to substantially reduce the number of frequencies initially needed to reconstruct the finite dimensional sources. By relating our analysis with results in [16], we also show that, with an additional requirement on the frequency set, reconstructions will be robust towards noise in the measurements. In Section 5, we conduct a set of numerical experiments to validate our findings.

## 2 The finite-dimensional source space and the singular value expansion

We now introduce and relate two important bases for $L^2(D_0)$, the Fourier-Bessel basis and the singular value expansion. This will motivate the construction of the finite-dimensional source space.

First, we describe the Fourier-Bessel (FB) basis in $L^2(D_0)$. This basis appears naturally when solving PDE problems in a disc geometry using separation of variables [23]. For $f \in L^2(D_0)$ we can expand

$$f(r,\theta) = \sum_{m\in\mathbb{Z}}\sum_{n\in\mathbb{N}} \hat{f}_{m,n}\varphi_{m,n}(r,\theta), \quad 0 \leq r \leq R_0, 0 \leq \theta \leq 2\pi, \tag{4}$$

where the angular-radial basis functions are given by

$$\varphi_{m,n}(r,\theta) = \left(\sqrt{\pi J_{m+1}^2(j_{m,n})}R_0\right)^{-1} e^{im\theta} J_m(k_{m,n}r). \tag{5}$$

$J_m(x)$ is the Bessel function of the first kind and of order $m$, and $j_{m,n}$ denotes the $n$'th, real and positive zero of $J_m(x)$, i.e., $J_m(j_{m,n}) = 0, \quad \forall m \in \mathbb{N}_0, n \in \mathbb{N}$. Here, $\mathbb{N}_0 = \mathbb{N} \cup \{0\}$. The numbers $k_{m,n}$ are the zeros $j_{m,n}$ scaled with respect to $R_0$, $k_{m,n} = j_{m,n}/R_0$. The equality in (4) is understood in the $L^2$-sense, and the functions $\{\varphi_{m,n}\}$ constitute an orthonormal basis for $L^2(D_0)$. The FB-coefficients $\hat{f}_{m,n}$ are computed with the standard inner product on $L^2(D_0)$, here in polar coordinates,

$$\hat{f}_{m,n} = (f,\varphi_{m,n})_{L^2(D_0)} = \int_0^{R_0}\int_0^{2\pi} f(r,\theta)\overline{\varphi_{m,n}(r,\theta)}rdrd\theta. \tag{6}$$

Since $J_{-m}(x) = (-1)^m J_m(x)$ for $m \in \mathbb{N}_0$, we have that[1] $\varphi_{-m,n} = (-1)^m\overline{\varphi_{m,n}}$. For real-valued $f \in L^2(D_0)$, it follows that

$$\overline{\hat{f}_{m,n}} = \overline{(f,\varphi_{m,n})}_{L^2(D_0)} = (-1)^m(f,\varphi_{-m,n})_{L^2(D_0)} = (-1)^m \hat{f}_{-m,n}. \tag{7}$$

Another consequence is that the scaled Bessel zero $k_{m,n}$ appears in both $\varphi_{m,n}$ and $\varphi_{-m,n}$.

The main idea in this paper relies on reconstructing source functions in a finite-dimensional subspace spanned by FB-functions. For $(M,N) \in \mathbb{N}_0 \times \mathbb{N}$ and $D_0$, we define the $(2M+1)N$-dimensional subspace

$$S_{M,N} = \text{span}\{\varphi_{m,n}, \ m=-M,\ldots,M, \ n=1,\ldots,N\}. \tag{8}$$

The projection $P_{M,N}$ of $s \in L^2(D_0)$ onto $S_{M,N}$ is given by

$$P_{M,N}s = \sum_{m=-M}^{M}\sum_{n=1}^{N} \hat{s}_{m,n}\varphi_{m,n}, \quad \text{where } \hat{s}_{m,n} = (s,\varphi_{m,n})_{L^2(D_0)}. \tag{9}$$

---

[1] The normalization constant is the same for indices $m$ and $-m$ since $\int J_{-m}(k_{m,n}r)^2 rdr = \int ((-1)^m J_m(k_{m,n}r))^2 rdr = \int J_m(k_{m,n}r)^2 rdr$.

Next, we introduce the singular value expansion (SVE) of the forward operator $F_k$ as given in [6, 16]. We write

$$F_k = \sigma_0^k(\cdot, \psi_0^k)_{L^2(D_0)}\phi_0^k + \sum_{n\in\mathbb{N}} \sigma_n^k\left((\cdot, \psi_n^k)_{L^2(D_0)}\phi_n^k + (\cdot, \psi_{-n}^k)_{L^2(D_0)}\phi_{-n}^k\right), \tag{10}$$

and the singular values and vectors, here given in polar coordinates, take the form

$$\sigma_n^k = \sqrt{2R}\pi R_0 |H_n^{(1)}(\kappa)| A_n(\kappa_0), \tag{11}$$

$$\psi_n^k(r,\theta) = (\sqrt{\pi}R_0 A_n(\kappa_0))^{-1} J_n(kr)e^{in\theta}, \tag{12}$$

$$\phi_n^k(\theta) = (2\pi R)^{-\frac{1}{2}} e^{i\arg H_n^{(1)}(\kappa)} e^{in\theta}, \tag{13}$$

where $n \in \mathbb{Z}$, $\kappa = kR$, $\kappa_0 = kR_0$, $H_n^{(1)}(z)$ is the Hankel function of the first kind and order $n$, and $A_n(\kappa) = \sqrt{J_n^2(\kappa) - J_{n-1}(\kappa)J_{n+1}(\kappa)}$. The singular vectors (12) and (13) constitute an orthonormal basis for $\ker(F_k)^T$ and $\text{range}(F_k)$, respectively.

Consider now the particular frequency $k = k_{m,n}$. By using the standard recursion formula for Bessel functions [25] and the fact that $\kappa_0 = j_{m,n}$, we get that

$$A_{\pm m}(\kappa_0) = \sqrt{-J_{m+1}(\kappa_0)J_{m-1}(\kappa_0)} = \sqrt{-J_{m+1}(\kappa_0)\left(\frac{2m}{\kappa_0}J_m(\kappa_0) - J_{m+1}(\kappa_0)\right)} = \sqrt{J_{m+1}^2(j_{m,n})}, \tag{14}$$

and hence

$$\psi_{\pm m}^{k_{m,n}} = (\sqrt{\pi}R_0 A_{\pm m}(\kappa_0))^{-1} J_{\pm m}(k_{m,n}r)e^{\pm im\theta} = \left(\sqrt{\pi J_{m+1}^2(j_{m,n})}R_0\right)^{-1} J_{\pm m}(k_{m,n}r)e^{\pm im\theta} = \varphi_{\pm m,n}. \tag{15}$$

Thus, at the particular frequencies $k_{m,n}$, the singular vectors $\psi_m^{k_{m,n}}$ coincide with the FB basis functions $\varphi_{m,n}$. We now note the following: given $s = \sum_{m,n} \hat{s}_{m,n}\varphi_{m,n} \in L^2(D_0)$ and the measurement $U_{k_{m,n}} = F_{k_{m,n}}s$ at the frequency $k_{m,n}$, we can invert the SVE to find the two FB-coefficients $\hat{s}_{\pm m,n}$, since

$$\frac{(U_{k_{m,n}}, \phi_{\pm m}^{k_{m,n}})_{L^2(\partial D)}}{\sigma_{|m|}^{k_{m,n}}} = (s, \psi_{\pm m}^{k_{m,n}})_{L^2(D_0)} = (s, \varphi_{\pm m,n})_{L^2(D_0)} = \hat{s}_{\pm m,n}. \tag{16}$$

Hence, if we have measurements at all frequencies $\{k_{m,n}\}$ we could use (16) to reconstruct any $s \in L^2(D_0)$. This is similar to what was proposed in [11].

We now consider a finite-dimensional setup. Define the frequency set consisting of the scaled Bessel zeros,

$$Q_{M,N} = \{k_{m,n}\}_{m=0,n=1}^{M,N}.$$

If we have measurements at all frequencies in the set $Q_{M,N}$, we could similarly reconstruct any source $s \in S_{M,N}$, or reconstruct the projection $P_{M,N}s$ of a general source $s \in L^2$.

This method would require one measurement $U_{k_{m,n}}$ for each reconstructed coefficient $\hat{s}_{m,n}$; but note that for each $(m,n)$ only one SVE-coefficient of $U_{k_{m,n}}$ is used. This gives rise to new questions of economical nature: From a measurement at some frequency $k$, how many, and which, FB-coefficients can be stably reconstructed? And given a finite-dimensional source space $S_{M,N}$, can we find a possibly small set of frequencies $Q_s = \{k_j\}$ such that we can guarantee numerically stable reconstruction of any source in $S_{M,N}$ from the subset of measurements $\{U_{k_j}\}_{k_j \in Q_s}$? We will address the latter question first, and the answer to the first question will follow.

## 3 Reconstruction of FB-coefficients from multi-frequency measurements

In this section, we obtain conditions on the set of frequencies needed to reconstruct a source in $S_{M,N}$. For multi-frequency measurements, we obtain this by constructing a change of basis matrix $K$ from the FB-basis to the SVE-basis. We then show that requiring that $K$ be invertible leads to certain conditions on the frequency set. In the next section, we show how knowledge about the distribution of the zeros of Bessel functions and the conditions on the frequencies allows us to substantially reduce the frequency set and still be able to stably reconstruct the source.

### 3.1 Change of basis: From SVE- to FB-coefficients

Let $Q_r = \{\tilde{k}_{m,n}\}_{m=0,n=1}^{M,N}$ be a set of *not necessarily distinct* frequencies. The indexing resembles that of $Q_{M,N}$; this is to emphasise that the two are closely connected, and to ease the some of the preceding arguments. Next, let $\{F_{\tilde{k}}\}_{\tilde{k} \in Q_r}$ be the family of forward operators associated with $Q_r$. We now define the space $V_{M,N} \subset L^2(D_0)$ to be

$$V_{M,N} = \operatorname{span}\{\psi_{\pm m}^{\tilde{k}_{m,n}} : \tilde{k}_{m,n} \in Q_r, m = 0, 1, \ldots, M, n = 1, 2, \ldots, N\}. \tag{17}$$

For a source $s \in S_{M,N}$, we have the FB-coefficients $\hat{s}_{m,n}$, and we introduce the SVE-coefficients $\hat{u}_{i,j} = (s, \psi_i^{\tilde{k}_{|i|,j}})_{L^2(D_0)}$. By expanding $s$ in the FB-basis we get the linear relation between the two sets of coefficients

$$\hat{u}_{i,j} = (s, \psi_i^{\tilde{k}_{|i|,j}})_{L^2(D_0)} = \sum_{m=-M}^{M} \sum_{n=1}^{N} \hat{s}_{m,n} (\varphi_{m,n}, \psi_i^{\tilde{k}_{|i|,j}})_{L^2(D_0)} = \sum_{n=1}^{N} \hat{s}_{i,n} (\varphi_{i,n}, \psi_i^{\tilde{k}_{|i|,j}})_{L^2(D_0)}, \tag{18}$$

where the last equality follows from the partial orthogonality relation

$$(\varphi_{m,n}, \psi_i^{\tilde{k}_{|i|,j}})_{L^2(D_0)} = \delta_{m,i} C(\tilde{k}_{|i|,j}, m, n, i)$$

where $C(\tilde{k}_{|i|,j}, m, n, i)$ is a complex number and $\delta_{m,i}$ is the Kronecker delta. We also note that

$$(\varphi_{-m,n}, \psi_{-m}^{\tilde{k}_{m,n}})_{L^2(D_0)} = (\varphi_{m,n}, \psi_m^{\tilde{k}_{m,n}})_{L^2(D_0)}, \tag{19}$$

due to a similar argument as in (7). Equation (18) constitutes a linear relation between the coefficients in the two bases. Next, define the vectors

$$\hat{S} = (\hat{s}_{1,1}, \hat{s}_{1,2}, \ldots, \hat{s}_{1,N}, \hat{s}_{2,1}, \ldots, \hat{s}_{M,N}, \hat{s}_{0,1}, \ldots, \hat{s}_{0,N}, \hat{s}_{-1,1}, \ldots, \hat{s}_{-M,N})^T,$$
$$\hat{U} = (\hat{u}_{1,1}, \hat{u}_{1,2}, \ldots, \hat{u}_{1,N}, \hat{u}_{2,1}, \ldots, \hat{u}_{M,N}, \hat{u}_{0,1}, \ldots, \hat{u}_{0,N}, \hat{u}_{-1,1}, \ldots, \hat{u}_{-M,N})^T,$$

and the matrix block-diagonal matrix $K \in \mathbb{R}^{(2M+1)N \times (2M+1)N}$ by

$$K = \begin{bmatrix} K^+ & 0 & 0 \\ 0 & K^0 & 0 \\ 0 & 0 & K^- \end{bmatrix}, \tag{20}$$

with $K^\pm \in \mathbb{R}^{MN \times MN}$ and $K^0 \in \mathbb{R}^{N \times N}$ defined as the matrices

$$K^+ = \begin{bmatrix} K_1^+ & 0 & \cdots & 0 \\ 0 & K_2^+ & \cdots & 0 \\ \vdots & \vdots & \ddots & \vdots \\ 0 & 0 & \cdots & K_M^+ \end{bmatrix}, \quad K^- = \begin{bmatrix} K_1^- & 0 & \cdots & 0 \\ 0 & K_2^- & \cdots & 0 \\ \vdots & \vdots & \ddots & \vdots \\ 0 & 0 & \cdots & K_M^- \end{bmatrix}, \tag{21}$$

$$K_m^+ = \{(\varphi_{m,n}, \psi_m^{\tilde{k}_{m,i}})_{L^2(D_0)}\}_{i=1,n=1}^{N,N} \quad \text{for } m = 1, 2, ...M, \tag{22}$$

$$K_m^- = \{(\varphi_{-m,n}, \psi_{-m}^{\tilde{k}_{m,i}})_{L^2(D_0)}\}_{i=1,n=1}^{N,N} \quad \text{for } m = 1, 2, ...M, \tag{23}$$

$$K^0 = \{(\varphi_{0,n}, \psi_0^{\tilde{k}_{0,i}})_{L^2(D_0)}\}_{i=1,n=1}^{N,N}, \tag{24}$$

where all $\tilde{k}_{m,i} \in Q_r$. The linear system (18) can now be written as $\hat{U} = K\hat{S}$. It is clear that $K$ is a mapping from the space of FB-coefficients for $S_{M,N}$ to the space of SVE-coefficients for the space $V_{M,N}$. The SVE-coefficients are related to the measurements $\{U_{k_j}\}$ through the relation

$$\hat{u}_{i,j} = \frac{(U_{\tilde{k}_{|i|,j}}, \phi_i^{\tilde{k}_{|i|,j}})_{L^2(\partial D)}}{\sigma_{|i|}^{\tilde{k}_{|i|,j}}}. \tag{25}$$

We see that reconstruction of $s$ amounts to reconstructing the vector $\hat{S}$ from $\hat{U}$, and hence depends on the invertibility of $K$. The properties of $K$ depend on the choice of frequency set $Q_r$. From (16), it is evident that if $Q_r = Q_{M,N}$, we have that $K = I$ since

$$\hat{u}_{i,j} = (s, \psi_i^{\tilde{k}_{|i|,j}})_{L^2(D_0)} = \sum_{n=1}^{N} \hat{s}_{i,n} (\varphi_{i,n}, \psi_i^{\tilde{k}_{|i|,j}})_{L^2(D_0)} = \hat{s}_{i,j}.$$

In the case that $Q_r \neq Q_{M,N}$, then $K \neq I$, but might still be invertible, and hence all the FB-coefficients should be recoverable from the SVE-coefficients. We therefore proceed to analyze the interaction between $Q_r$ and the invertibility of $K$.

### 3.2 Invertibility of $K$

This section contains Theorem 1, which is our main result. It says that if the frequencies in $Q_r$ are not too far away from those in $Q_{M,N}$, then $K$ is invertible. Lemma 1 gives a computable bound for this distance, and is required in the proof of Theorem 1. We will need the functions

$$\mu_m = \begin{cases} \frac{R_0}{\pi} & \text{for } m \geq 1, \\ \frac{R_0}{j_{0,2} - j_{0,1}} & \text{for } m = 0, \end{cases} \quad \text{and} \quad \chi_{j,i} = \begin{cases} 1 & \text{for } i \neq j, \\ 0 & \text{for } i = j. \end{cases} \tag{26}$$

**Lemma 1.** *For $R_0 > 0$ and $(M, N) \in \mathbb{N}_0 \times \mathbb{N}$, let $Q_{M,N} = \{j_{m,n}/R_0\}_{m=0,n=1}^{M,N}$. Then, for any indices $0 \leq m \leq M$ and $1 \leq i \leq N$, there is a largest positive number $\Delta k_{m,i}$ such that*

$$k_{m,i} \pm \Delta k_{m,i} \in (k_{m,i-1}, k_{m,i+1})$$

*and that the inequality*

$$\left| \frac{k_{m,i}}{(k_{m,i} \pm \delta)^2 - k_{m,i}^2} \right| \geq \chi_{1,i} \left( \frac{\mu_m}{2} \log \left( \frac{(k_{m,i} \pm \delta)^2 - k_{m,1}^2}{(k_{m,i} \pm \delta)^2 - k_{m,i-1}^2} \right) + \frac{k_{m,i-1}}{(k_{m,i} \pm \delta)^2 - k_{m,i-1}^2} \right)$$
$$+ \chi_{N,i} \left( \frac{\mu_m}{2} \log \left( \frac{k_{m,N}^2 - (k_{m,i} \pm \delta)^2}{k_{m,i+1}^2 - (k_{m,i} \pm \delta)^2} \right) + \frac{k_{m,i+1}}{k_{m,i+1}^2 - (k_{m,i} \pm \delta)^2} \right) \quad (27)$$

*holds for all $\delta \in (0, \Delta k_{m,i})$.*

*Furthermore, let $\Delta k = \min_{m,i} \Delta k_{m,i}$. Then the inequality*

$$\left| \frac{k_{m,i}}{(k_{m,i} \pm \delta)^2 - k_{m,i}^2} \right| > \sum_{n=1, n \neq i}^{N} \left| \frac{k_{m,n}}{(k_{m,i} \pm \delta)^2 - k_{m,n}^2} \right| \quad (28)$$

*holds for all $0 < \delta < \Delta k$ and all indices $0 \leq m \leq M$ and $1 \leq i \leq N$.*

*Proof.* For $R_0 > 0$ and $(M, N) \in \mathbb{N}_0 \times \mathbb{N}$, let $Q_{M,N} = \{k_{m,n}\}_{m=0,n=1}^{M,N}$, where $k_{m,n} = j_{m,n}/R_0$. For fixed $0 \leq m \leq M$, the sequence $\{k_{m,n}\}_{n=1}^{N}$ is strictly increasing and positive (see e.g. Proposition 1 in Section 4.1). We fix $m$ and $i$, and look at the inequality

$$\left| \frac{k_{m,i}}{x^2 - k_{m,i}^2} \right| > \sum_{n=1, n \neq i}^{N} \left| \frac{k_{m,n}}{x^2 - k_{m,n}^2} \right| \quad (29)$$

with the assumption $k_{m,i-1} < x < k_{m,i+1}$. We estimate the sum of the right-hand side. Define the functions

$$S^-(x, k) = \frac{k}{x^2 - k^2} \quad \text{for } 0 < k < x, \quad \text{and} \quad S^+(x, k) = \frac{k}{k^2 - x^2} \quad \text{for } x < k,$$

and write

$$\sum_{n=1, n \neq i}^{N} \left| \frac{k_{m,n}}{x^2 - k_{m,n}^2} \right| = \sum_{n=1}^{i-1} S^-(x, k_{m,n}) + \sum_{n=i+1}^{N} S^+(x, k_{m,n}).$$

Since $S^-(x, k_{m,n}) \leq S^-(x, s)$ for $k_{m,n} \leq s < x$, and $S^+(x, k_{m,n}) \leq S^+(x, s)$ for $x \leq s < k_{m,n}$, we obtain the inequalities

$$S^-(x, k_{m,n}) < \mu_m \int_{k_{m,n}}^{k_{m,n+1}} \frac{s}{x^2 - s^2} ds$$
$$= \frac{\mu_m}{2} (\log(x^2 - k_{m,n}^2) - \log(x^2 - k_{m,n+1}^2)), \quad \text{for } x < k_{m,i-1},$$

$$S^+(x, k_{m,n}) < \mu_m \int_{k_{m,n-1}}^{k_{m,n}} \frac{s}{s^2 - x^2} ds$$
$$= \frac{\mu_m}{2} (\log(k_{m,n}^2 - x^2) - \log(k_{m,n-1}^2 - x^2)), \quad \text{for } x > k_{i+1}.$$

The multiplication by $\mu_m$, defined in (26), corrects for the integration being done over an interval of length $k_{m,i+1} - k_{m,i}$, since $k_{m,i+1} - k_{m,i} \geq \frac{\pi}{R_0}$ for $m, i \geq 1$, and $k_{0,i+1} - k_{0,i} \geq \frac{j_{0,2} - j_{0,1}}{R_0}$ for $i \geq 1$. These properties are described in Section 4, where we go into more detail about the properties of the zeros of Bessel functions. The estimates are quite sharply away from the diagonal, except for $S^-(x, k_{m,n-1})$ and $S^+(x, k_{m,n+1})$, so we include these terms directly in the sums and get

$$\sum_{n=1}^{i-1} S^-(x, k_{m,n}) < \frac{\mu_m}{2}(\log(x^2 - k_{m,1}^2) - \log(x^2 - k_{m,i-1}^2)) + S^-(x, k_{m,i-1})$$

$$= \frac{\mu_m}{2} \log\left(\frac{x^2 - k_{m,1}^2}{x^2 - k_{m,i-1}^2}\right) + \frac{k_{m,i-1}}{x^2 - k_{m,i-1}^2}$$

$$\sum_{n=i+1}^{N} S^+(x, k_{m,n}) < \frac{\mu_m}{2}(\log(k_{m,N}^2 - x^2) - \log(k_{m,i+1}^2 - x^2)) + S^+(x, k_{m,i+1})$$

$$= \frac{\mu_m}{2} \log\left(\frac{k_{m,N}^2 - x^2}{k_{m,i+1}^2 - x^2}\right) + \frac{k_{m,i+1}}{k_{m,i+1}^2 - x^2}.$$

To have the estimate valid for the cases $i = 1$ and $i = N$, we need the function $\chi_{i,j}$ defined in (26). We then have that

$$\sum_{n=1, n \neq i}^{N} \left|\frac{k_{m,n}}{x^2 - k_{m,n}^2}\right| < \chi_{1,i} \left(\frac{\mu_m}{2} \log\left(\frac{x^2 - k_{m,1}^2}{x^2 - k_{m,i-1}^2}\right) + \frac{k_{m,i-1}}{x^2 - k_{m,i-1}^2}\right)$$
$$+ \chi_{N,i} \left(\frac{\mu_m}{2} \log\left(\frac{k_{m,N}^2 - x^2}{k_{m,i+1}^2 - x^2}\right) + \frac{k_{m,i+1}}{k_{m,i+1}^2 - x^2}\right). \quad (30)$$

We rewrite $x$ as a perturbation of $k_{m,i} \in Q_{M,N}$ on the form $x = k_{m,i} \pm \delta$, with $\delta > 0$ and such that $x \in (k_{m,i-1}, k_{m,i+1})$. Let $\Delta k_{m,i}$ be the largest number such that

$$\left|\frac{k_{m,i}}{(k_{m,i} \pm \delta)^2 - k_{m,i}^2}\right| \geq \chi_{1,i} \left(\frac{\mu_m}{2} \log\left(\frac{(k_{m,i} \pm \delta)^2 - k_{m,1}^2}{(k_{m,i} \pm \delta)^2 - k_{m,i-1}^2}\right) + \frac{k_{m,i-1}}{(k_{m,i} \pm \delta)^2 - k_{m,i-1}^2}\right)$$
$$+ \chi_{N,i} \left(\frac{\mu_m}{2} \log\left(\frac{k_{m,N}^2 - (k_{m,i} \pm \delta)^2}{k_{m,i+1}^2 - (k_{m,i} \pm \delta)^2}\right) + \frac{k_{m,i+1}}{k_{m,i+1}^2 - (k_{m,i} \pm \delta)^2}\right), \quad (31)$$

holds for all $\delta \in (0, \Delta \kappa_{m,i})$. Such a $\Delta \kappa_{m,i}$ exists, since for fixed $m$ and $i$ and a closed interval $[k_{m,i} - \delta, k_{m,i} + \delta] \subset (k_{m,i-1}, k_{m,i+1})$ the estimate on the right-hand side of (30) is continuous and bounded, and hence attains a maximum $M$. On the other hand, the function $|k_{m,i}/((k_{m,i} \pm \delta)^2 - k_{m,i}^2)|$ is unbounded from both sides as $\delta \to 0$, and we can find a $\delta$ such that $\min\{|k_{m,i}/((k_{m,i} \pm \delta)^2 - k_{m,i}^2)|\} \geq M$.

It follows from (30) that

$$\left|\frac{k_{m,i}}{(k_{m,i} \pm \delta)^2 - k_{m,i}^2}\right| > \sum_{n=1, n \neq i}^{N} \left|\frac{k_{m,n}}{(k_{m,i} \pm \delta)^2 - k_{m,n}^2}\right| \quad (32)$$

holds for all $\delta \in (0, \Delta k_{m,i})$, and that with $\Delta k = \min_{m,i} \Delta k_{m,i}$, it holds for all indices $0 \leq m \leq M$ and $1 \leq i \leq N$. $\square$

We are now in the position to prove our main result.

**Theorem 1.** *Consider an ISP geometry specified by $R$ and $R_0$, and a finite-dimensional source space $S_{M,N}$, and let $\Delta k > 0$ be as specified in Lemma 1. If then, for each $\tilde{k}_{m,i} \in Q_r$ and corresponding $k_{m,i} \in Q_{M,N}$, we have that $|\tilde{k}_{m,i} - k_{m,i}| < \Delta k$, the matrix $K$ is invertible, and any $s \in S_{M,N}$ is uniquely recoverable from measurements at the frequencies $\tilde{k} \in Q_r$.*

*Proof.* We will show that under the assumption in Theorem 1, $K$ is strictly diagonally dominant. Recall that a matrix $A \in \mathbb{R}^{m \times m}$ is called strictly diagonally dominant if

$$|a_{i,i}| > \sum_{j=1, j \neq i}^{m} |a_{i,j}| \quad \text{for } i = 1, 2, ...m. \tag{33}$$

It is well-known, see e.g. [24], that strictly diagonally dominant matrices are invertible.

We start by computing the entries of matrix $K$, as given in (22)-(24). For the matrix-block $K_m^+$, the entries of the $i$'th row are given as

$$(\varphi_{m,n}, \psi_m^{\tilde{k}_{m,i}})_{L^2(D_0)}, \quad \text{for } n = 1, 2, ..., N.$$

Using (5) and (10), we compute

$$(\varphi_{m,n}, \psi_m^{\tilde{k}_{m,i}})_{L^2(D_0)} = \int_0^{R_0} \int_0^{2\pi} \psi_m^{\tilde{k}_{m,i}} \overline{\varphi_{m,n}} r \mathrm{d}\theta \mathrm{d}r$$

$$= \frac{1}{\pi R_0^2 \sqrt{J_{m+1}^2(j_{m,n})} A_m(R_0 \tilde{k}_{m,i})} \int_0^{2\pi} e^{im\theta} e^{-im\theta} \int_0^{R_0} J_m(\tilde{k}_{m,i} r) J_m(k_{m,n} r) r \mathrm{d}\theta \mathrm{d}r$$

$$= \frac{2}{R_0^2 \sqrt{J_{m+1}^2(j_{m,n})} A_m(R_0 \tilde{k}_{m,i})} \int_0^{R_0} J_m(\tilde{k}_{m,i} r) J_m(k_{m,n} r) r \mathrm{d}r.$$

The partial orthogonality appearing in (18) is apparent from the second line. The evaluation of the last integral can be found in, e.g., [25], and with some further manipulation we get that

$$(\varphi_{m,n}, \psi_m^{\tilde{k}_{m,i}})_{L^2(D_0)} = -\mathrm{sign}(J_{m+1}(j_{m,n})) \frac{2 J_m(\tilde{k}_{m,i} R_0)}{R_0 A_m(\tilde{k}_{m,i} R_0)} \frac{k_{m,n}}{(\tilde{k}_{m,i}^2 - k_{m,n}^2)}. \tag{34}$$

Due to the block diagonal structure of $K$, and since (19) implies that $K_m^+ = K_m^-$, it is sufficient to check that each block $K_m^+$ is diagonally dominant. On the $i$'th row of $K_m^+$, the inequality in (33) takes the form

$$|(\varphi_{m,i}, \psi_m^{\tilde{k}_{m,i}})_{L^2(D_0)}| > \sum_{n=1, n \neq i}^{N} |(\varphi_{m,n}, \psi_m^{\tilde{k}_{m,i}})_{L^2(D_0)}|, \tag{35}$$

or, in terms of (34),

$$\left| \frac{2 J_m(\tilde{k}_{m,i} R_0)}{R_0 A_m(\tilde{k}_{m,i} R_0)} \frac{k_{m,i}}{(\tilde{k}_{m,i}^2 - k_{m,i}^2)} \right| > \sum_{n=1, n \neq i}^{N} \left| \frac{2 J_m(\tilde{k}_{m,i} R_0)}{R_0 A_m(\tilde{k}_{m,i} R_0)} \frac{k_{m,n}}{(\tilde{k}_{m,i}^2 - k_{m,n}^2)} \right|. \tag{36}$$

We note that, due to the relation in (15), we have that when $\tilde{k}_{m,i} \to k_{m,i}$, the left-hand side of (36) equals 1, while the right-hand side equals 0. Hence (36) holds. When $\tilde{k}_{m,i} \neq k_{m,i}$, the factor $(2J_m(\tilde{k}_{m,i}R_0))/(R_0A_m(\tilde{k}_{m,i}R_0))$ enters as a constant on both sides, and can therefore be neglected. Hence we need only look at the inequality

$$\left|\frac{k_{m,i}}{\tilde{k}_{m,i}^2 - k_{m,i}^2}\right| > \sum_{n=1,n\neq i}^{N} \left|\frac{k_{m,n}}{\tilde{k}_{m,i}^2 - k_{m,n}^2}\right|. \tag{37}$$

Recognizing the above expression as equation (28) in Lemma 1, we see that if $|\tilde{k}_{m,i} - k_{m,i}| < \Delta k$, with $\Delta k$ as specified in Lemma 1, then (37) must hold. It follows that if $|\tilde{k}_{m,i} - k_{m,i}| < \Delta k$ holds for all $0 \leq m \leq M$ and $1 \leq i \leq N$, then (37) holds for each row in each $K_m^\pm$, and that $K$ is strictly diagonally dominant and hence invertible. $\square$

In Figure 1 we plot both sides of (27) and the right-hand side of (28), for different choices of $m$ and $i$ and fixed $N$. We see that the interval where the inequality holds is quite insensitive to changes in indices. Also, numerical investigation indicates that when $M, N \leq 150$, it is sufficient to check the inequality only for $m = 0$ and $n = 1$, i.e., the first row in $K$, since the smallest admissible $\Delta k$ occurs there (We did not try with $M, N > 150$). We suspect this is related to the convexity/concavity of the zeros, and a proof of this would definitely be interesting.

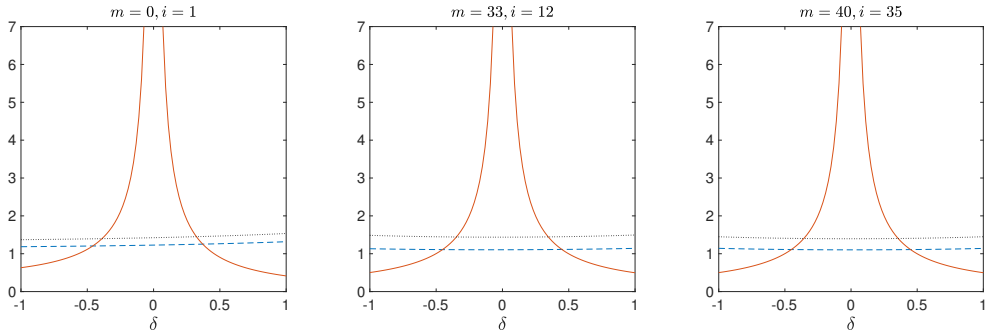

Figure 1: The red and the (dotted) black lines are the left- and right-hand sides of (31), respectively. The dashed blue line is the estimated off-diagonal sum $\sum_{n=1,n\neq i}^{N} \left|k_{m,n}/(\tilde{k}_{m,i}^2 - k_{m,n}^2)\right|$. We use $\delta \in [-1, 1]$, and set $R_0 = 1$ and $N = 50$.

Lemma 1 gives a bound $\Delta k$ on admissible perturbations of frequencies in $Q_{M,N}$, sufficient to guarantee invertibility of $K$. With the set $Q_{M,N}$ available it is straightforward to compute $\Delta k$. In Section 4, we analyze the distribution of the zeros of Bessel functions, and show how one can choose a set of frequencies $Q_s$ that are admissible perturbations of the frequencies $Q_{M,N}$, but such that $Q_s$ is substantially smaller than $Q_{M,N}$. The set $Q_r$ in Theorem 1 is then constructed using the frequencies in $Q_s$.

### 3.3 Reconstruction of sources not in $S_{M,N}$

When we take the source to be any function $s \in L^2(D_0)$, we introduce an error in the reconstruction of the projection $s$ on $S_{M,N}$. Let $U_k = F_k s$. The error occurs in the SVE-coefficients; Let $\tilde{u}_{i,j}$ be a coefficient computed from a measurement of source in $L^2$. From (18) we then have

$$\tilde{u}_{i,j} = \frac{(U_k, \phi_m^k)_{L^2(\partial D_0)}}{\sigma_{|i|}^k} = (s, \psi_i^k)_{L^2(D_0)} = \sum_{n=1}^{\infty} \hat{s}_{i,n}(\varphi_{i,n}, \psi_i^k)_{L^2(D_0)}. \tag{38}$$

Splitting the last sum, we get

$$\tilde{u}_{i,j} = \sum_{n=1}^{N} \hat{s}_{i,n} (\varphi_{i,n}, \psi_i^k)_{L^2(D_0)} + \sum_{n=N+1}^{\infty} \hat{s}_{i,n} (\varphi_{i,n}, \psi_i^k)_{L^2(D_0)} = \hat{u}_{i,j} + E(i,N). \tag{39}$$

Hence, we get an additional error term $E(i,N)$ in our coefficients, and the reconstruction of the projection of $s$ on $S_{M,N}$ will not be exact. If the source is assumed to have certain regularity properties, one can estimate on the size of $E$, but in general one cannot. However, in numerical simulations, the effect of $E(i,N)$ on reconstructions appears to be small, even for discontinuous sources.

# 4 Invertibility and stability with a reduced number of frequencies

In this section, we investigate the consequences of the findings in Section 3 for a reduction of the measurement frequencies. By considering the distribution of the zeros of Bessel functions, we will see that a lot of the frequencies in $Q_{M,N}$ are redundant. We also analyze the stability and regularizing effect that comes with considering finite-dimensional reconstructions.

## 4.1 Density of the zeros of Bessel functions in $\mathbb{R}^+$

As a consequence of the findings in Section 3, we are interested in the distribution of the positive zeros of Bessel functions on the real line. We want to investigate frequencies that are in the neighborhood of such zeros, particularly frequencies that are in the neighborhood of multiple zeros. Given positive $\alpha$ and $\Delta j$, and the interval $I_\alpha(\Delta j) = (\alpha - \Delta j, \alpha + \Delta j)$, we want to know how many zeros $j_{m,n}$ are contained within $I_\alpha(\Delta j)$, and how this depends on $\alpha$ and $\Delta j$.

The literature on Bessel functions is vast, and we gather here a few properties of their zeros that proved useful for our purpose. They are found in [8], [10],[22] and [20], respectively.

**Proposition 1.**

1) $\lim_{n \to \infty} (j_{m,n+1} - j_{m,n}) = \pi$ for $m \in \mathbb{R}$.

2) $j_{m,n} > m + n\pi - \frac{\pi}{2} + \frac{1}{2}, \quad m \geq -1/2$.

3) If $|m| \neq \frac{1}{2}$, then $(j_{m,n+1} - j_{m,n} - \pi)(|m| - \frac{1}{2}) > 0$.

4) Interlacing property: $j_{m,1} < j_{m+1,1} < j_{m,2} < j_{m+1,2} < j_{m,3} < \cdots$

From 2) it is evident that the zeros are strictly increasing functions of both index and order, and that they are shifted upwards on the positive real axis with increasing order. Also, 1) and 3) show that the distribution of the zeros for a fixed order becomes, eventually, quite regular, as the difference between consecutive zeros tends towards $\pi$. The interlacing property 4) shows how zeros are spaced between each other as a function of order and index. This suggests that a zero $j_{m,n}$ may be close to some zero $j_{i,j}$, or $j_{m,n} \approx j_{i,j}$, when $m > i$ and $j > n$, or vice versa. In Figure 2, we have plotted the 7 first zeros of $J_m$, $m = 0, 1, 2, ..., 7$. The horizontal lines indicate an interval $I_\alpha(\Delta j)$. Observe that different zeros of $J_m$'s of different order lie within the strip spanned by $I_\alpha(\Delta k)$ in the $\mathbb{R}^+ \times \mathbb{N}_0$-plane. Figure 3 shows the number of $j_{n,l}$ such that $|\alpha - j_{n,l}| < 0.5$ for an increasing $\alpha$. The numerical results indicate that the density of zeros increases as one moves up the real line. In

fact, one can produce an approximate formula for the number of close zeros to a given point from the line argument. The reasoning follows from considering a horizontal interval centered at a point $\alpha$ in the $\mathbb{R}^+ \times \mathbb{N}_0$-plane, as in Figure 2.

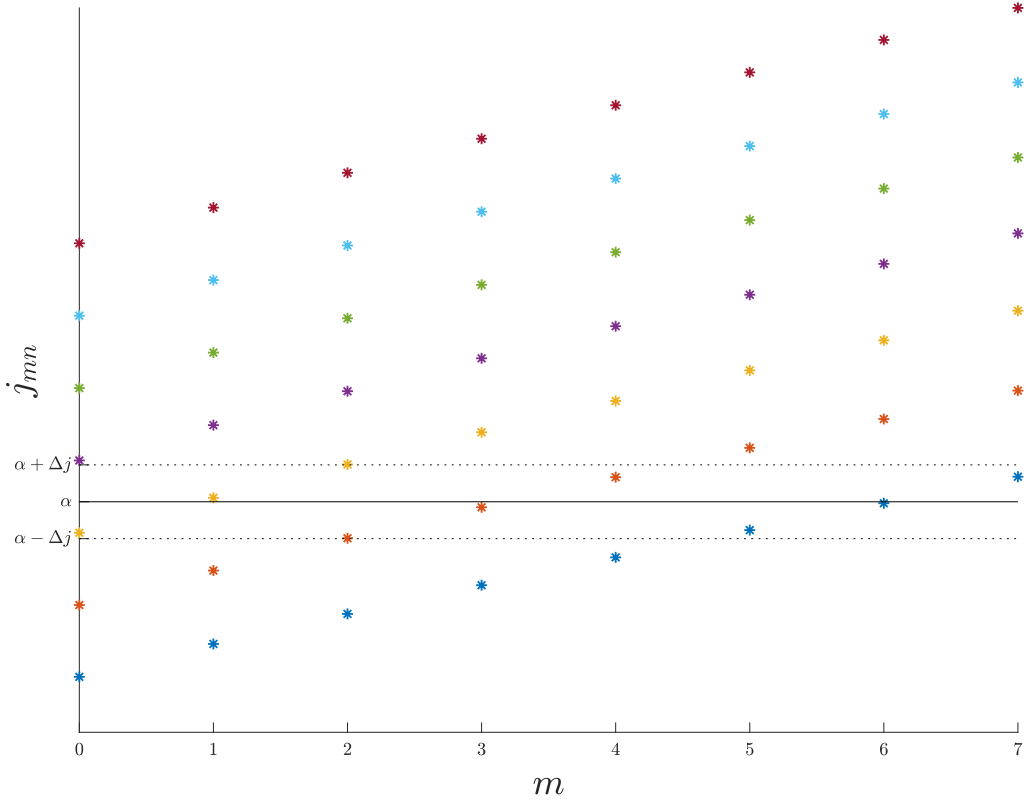

Figure 2: A plot of the 7 first zeros of $J_m$, $m = 0, 1, 2, ..., 7$. The dashed horizontal lines bound an interval $(\alpha - \Delta j, \alpha + \Delta j)$, with $\alpha = 10$ and $\Delta j = 1.6$.

A horizontal line passing through a point $\alpha > 0$ will intersect the intervals $(j_{m,1}, \infty)$ for all $m$ such that $\alpha > j_{m,1}$. Since the distance between consecutive zeros tends towards $\pi$, we can approximate the expected number of zeros of $J_m(x)$ in an interval $I_\alpha(\Delta j), \alpha \geq j_{m,1}$, to be $\frac{2\Delta j}{\pi}$, and note that this approximation will get better as $\alpha$ increases. We now estimate the number of intervals $(j_{m,1}, \infty)$ the horizontal line passing through $\alpha$ intersects. Since $j_{m,1}$ is an increasing function of $m$, this is equivalent to estimating the largest $M \in \mathbb{N}$ such that $\alpha \geq j_{M,1}$ holds. To this end we use the following upper bound on $j_{M,1}$ [25]:

$$j_{M,1} < \left(\frac{4}{3}(M+1)(M+5)\right)^{\frac{1}{2}}, M > 0.$$

Requiring $j_{M,1} < \alpha$, we now want the largest $M$ such that the following inequality holds:

$$\left(\frac{4}{3}(M+1)(M+5)\right)^{\frac{1}{2}} \leq \alpha. \tag{40}$$

Hence we get an estimate for $M$,

$$M = \left\lfloor \left( \frac{1}{2} \left(3\alpha^2 + 16\right)^{\frac{1}{2}} - 3 \right) \right\rfloor.$$

Combining $M$ and the approximation of the expected number of zeros, we have that the total number of zeros within the interval $[\alpha - \Delta j, \alpha + \Delta j]$, $N(\alpha, \Delta j)$, is approximately

$$N(\alpha, \Delta j) \approx M \frac{2\Delta j}{\pi} = \left\lfloor \left( \frac{1}{2} \left(3\alpha^2 + 16\right)^{\frac{1}{2}} - 3 \right) \right\rfloor \frac{2\Delta j}{\pi} \tag{41}$$

This approximation is also plotted in Figure 3. We see that $N(\alpha, \Delta j) = O(\alpha \Delta j)$. The formula implies that, for a fixed $\Delta j$, the number of zeros contained in an interval $I_\alpha(\Delta j)$ increases with increasing $\alpha$. This is in agreement with the numerical findings.

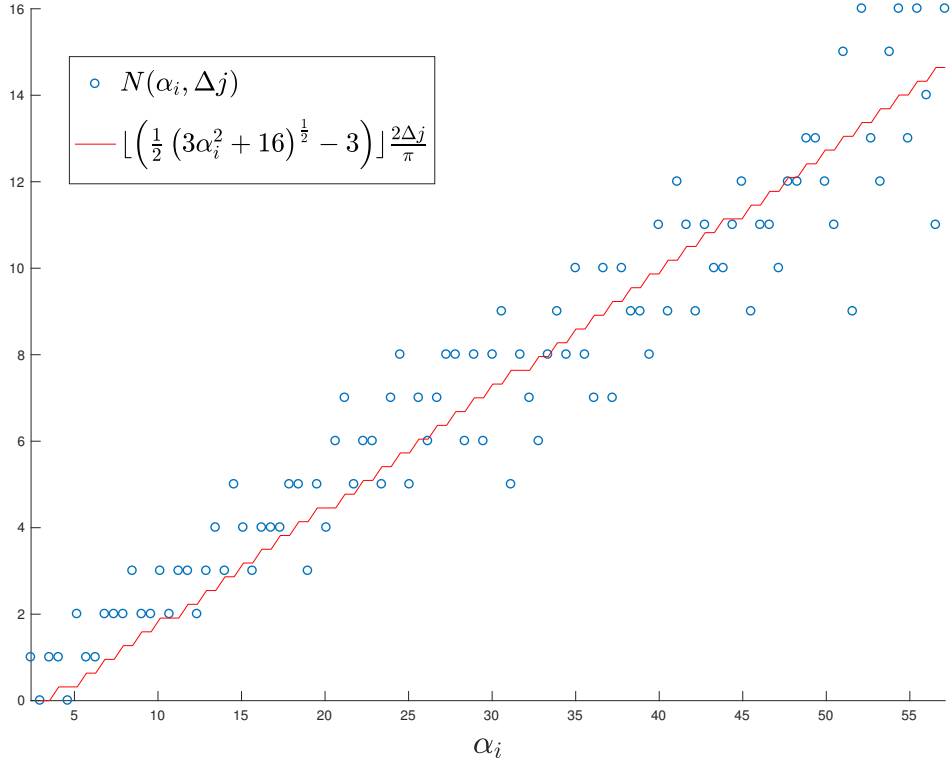

Figure 3: The plot shows $N(\alpha_i, 0.5)$ for $\alpha_i = j_{0,1} + j_{50,1}/(N+1)i = 2.4048 + 0.5769i, i = 0, 1, 2, ..., 99$. The red line is a plot of the approximative formula for $N(\alpha_i, (\Delta j))$ in (41).

To further substantiate the assertion that the zeros of the Bessel functions $J_m(x), m = 0, 1, 2...$ become more densely distributed as their numerical value increases, we include the following result: In [10, 15] it is shown that if we consider the set $\cup_{j=1}^\infty \{x_j\} = \cup_{m,n=1}^\infty \{j_{m,n}\}$, with $x_j$'s strictly increasing, then

$$\lim_{j \to \infty} x_{j+1} - x_j = 0.$$

Also, in [21], the authors did a numerical investigation of the distribution of the zeros of $J_m(x), m = 0, 1, 2, .....2000$ for $x \in [0, 3000]$, and computed the distances $d_i$ between the (sorted) zeros. They found that the maximum distance $d_i$ between two zeros was max $d_i = 1.426$, and that the mean and standard deviation where $\mu(d_i) = 0.0028$ and $\sigma(d_i) = 0.0062$, respectively.

Together, the above results imply that the answer to our initial question about the existence of frequencies that can are suitable perturbations of Bessel zeros, is positive. Furthermore, the results show that the number of zeros in a neighborhood of a given frequency increases as the frequency increases. Theorem 1 and the results in this section imply that we can reduce the number of different frequencies needed for $K$ to be invertible.

## 4.2 Reducing the frequency set

Once we have decided on a finite dimensional source space, $S_{M,N}$ we can use Theorem 1 to compute a perturbation level $\Delta k$ such that the SVE- to FB-coefficient mapping $K$ is invertible under perturbation of frequencies $Q_{M,N}$. Given $Q_{M,N}$, set $k^- = \min k_{m,n}$ and $k^+ = \max k_{m,n}$. We define $Q_s(\Delta k)$ to be *the smallest member of the family of sets*

$$\Big\{ k \in [k^-, k^+] : \text{ for every } k_{m,n} \in Q_{M,N} \quad \exists k \text{ such that } |k_{m,n} - k| < \Delta k \Big\}.$$

The set $Q_s$ is sometimes called a minimal subcover. Table 4.2 contains some numbers on the size of $Q_s(\Delta k)$ compared to $Q_{M,N}$ for different values of $\Delta k$ and $M$ and $N$. The increasing density of the Bessel zeros is evident, as the relative size of the reduced frequency set $Q_s$ is much smaller for the full frequency set $Q_{M,N}$ for larger $M$ and $N$.

| $\Delta k$ | | 0 | 0.25 | 0.5 | 0.75 | 1 | 1.25 |
|---|---|---|---|---|---|---|---|
| $M = N = 7$ | $\frac{|Q_s|}{|Q_{M,M}|} \cdot 100\%$ | 100% | 62.5% | 41.1% | 30.4% | 23.2% | 19.6% |
| $M = N = 7$ | $|Q_s|$ | 56 | 35 | 23 | 17 | 13 | 11 |
| $M = N = 50$ | $\frac{|Q_s|}{|Q_{M,M}|} \cdot 100\%$ | 100% | 15% | 8.2% | 5.6% | 4.2% | 3.5% |
| $M = N = 50$ | $|Q_s|$ | 2550 | 382 | 210 | 143 | 108 | 88 |

Table 1: The numerical value in the first row gives the reduction percentage for the reduced frequency set, and the second row contains the number of frequencies in the different $Q_s$'s. For dimensions $m = n = 7$ and $M = N = 50$ and $R_0 = 1$, one has $|Q_{m,n}| = 56$ and $|Q_{M,N}| = 2550$.

Once the reduced frequency set is computed, one can construct $K$ using the smaller number of perturbed frequencies in $Q_s$ instead of those in $Q_{M,N}$. We assemble the frequency set $Q_r$ of *not necessarily distinct frequencies* in Theorem 1 by substituting those $k_{m,n} \in Q_{M,N}$ by their closest $k_j \in Q_s$. Hence $Q_r$ will have several repeated entries.

## 4.3 Stability and regularization through the reduced frequency set

Next, we analyze the stability of the equation $K\hat{S} = \hat{U}$ with respect to additive noise in the measurement data. Note that the invertibility of $K$ is unaffected by noisy data. We assume a noisy measurement at frequency $k_j$ is of the form $U^\varepsilon_{k_j} = U_{k_j} + \varepsilon$, where $\varepsilon \in L^2(\partial D)$ is the noise and $U_{k_j}$ is the true measurement. Also, we assume the relative noise level is bounded by a $\delta > 0$, $\|\varepsilon\|_{L^2}/\|U_{k_j}\|_{L^2} \leq \delta$. To see how a noisy measurement influences our solution, we define the sets

$$I_{k_j} = \{(m,n) : k_{m,n} \in Q_{M,N} \text{ and } |k_{m,n} - k_j| < \Delta k\},$$
$$M_{k_j} = \{m : k_{m,n} \in Q_{M,N} \text{ and } |k_{m,n} - k_j| < \Delta k\}.$$

$I_{k_j}$ the set of indices of those $k_{m,n} \in Q_{M,N}$ that are perturbations of $k_j \in Q_s$, and $M_{k_j}$ contains only the order index $m$. Recalling equation (25) from the construction of the linear system $K\hat{S} = \hat{U}$, we see that $I_{k_j}$ contains the indices of the singular vectors $\phi_i^{k_j}$ and singular values $\sigma_i^{k_j}$ used to compute the SVE-coefficients $\hat{u}_{m,n}^\varepsilon$ from the measurement $U_{k_j}^\varepsilon$:

$$\hat{u}_{\pm m,n}^\varepsilon = \frac{(U_{k_j}^\varepsilon, \phi_{\pm m}^{k_j})_{L^2(\partial D)}}{\sigma_m^{k_j}}, \quad (m,n) \in I_{k_j}. \tag{42}$$

It is well-known that the stability of inverting a compact operator depends on the behavior of the singular value spectrum, and that small singular values amplify the (high-frequency) noise. Such is the case here:

$$|\hat{u}_{m,n}| = \left|\frac{(U_{k_j}^\varepsilon, \phi_m^{k_j})_{L^2}}{\sigma_m^{k_j}}\right| \leq \left|\frac{(U_{k_j}, \phi_m^{k_j})_{L^2}}{\sigma_m^{k_j}}\right| + \left|\frac{(\varepsilon, \phi_m^{k_j})_{L^2}}{\sigma_m^{k_j}}\right| \leq |\hat{u}_{m,n}| + \frac{\|\varepsilon\|_{L^2}}{\sigma_m^{k_j}}.$$

It is evident that when $\sigma_m \ll 1$, we get an amplification of the noise. Due to the fact [17] that $\lim_{m \to \infty} \sigma_m^{k_j} = 0$, it becomes important to study the pre-asymptotic behavior of the singular values. In [16], a characterization of the "bandwidth" of the operator $F_{k_j}$ is presented, where the bandwidth $\beta_{k_j}$ is defined as the singular value index whereafter the singular values of $F_{k_j}$ become strictly decreasing:

$$\beta_{k_j} = \mathrm{argmin}_{m \in \mathbb{N}_0}\{\sigma_{m+n}^{k_j} > \sigma_{m+n+1}^{k_j} \text{ for all } n \in \mathbb{N}_0\}.$$

Knowing the value of $\beta_{k_j}$, we can say whether noise in the SVE-coefficient $\hat{u}_{m,n}$ will be amplified. In [16], a tight lower bound for the bandwidth is proved, and a tight upper bound is conjectured. Both bounds are numerically validated. We state the result on the lower bound $M_-$:

The bandwidth $\beta_{k_j}$ of the forward operator $F_{k_j} : s \mapsto U_{k_j}$ associated with (1) satisfies $M_- \leq \beta_{k_j}$, with $M_-$ given by

$$M_- = \mathrm{argmin}_{m \in \mathbb{N}_0}\{j_{m,1} \geq k_j R_0\}.$$

Hence, to avoid the amplification of noise from a measurement $U_{k_j}^\varepsilon$, we only want to extract the following SVE-coefficients:

$$\hat{u}_{\pm m,n}^\varepsilon = \frac{(U_{k_j}^\varepsilon, \phi_{\pm m}^{k_j})_{L^2(\partial D)}}{\sigma_m^{k_j}}, \quad m \in M_{k_j}, (m,n) \in I_{k_j} \text{ and } m \leq \beta_{k_j}.$$

We now ask if we can avoid $m \in M_{k_j}$ such that $m > \beta_{k_j}$. The answer is given by the following proposition:

**Theorem 2.** *If $\Delta k \leq 1/R_0$, then all indices $m$ in $M_{k_j}$ are within the bandwidth of the forward operator $F_{k_j}$, that is, $\max_{M_{k_j}} m \leq \beta_{k_j}$.*

*Proof.* We first need to show that $j_{m+1,1} - j_{m,1} > 1$ for $m \in \mathbb{N}_0$. From [14], we have that $\frac{dj_{m,n}}{dm} > 1$ whenever $m > -1$ and $n = 1, 2, 3, \ldots$. Hence it follows by the fundamental theorem of calculus that

$$j_{m+1,1} - j_{m,1} = \int_m^{m+1} \frac{dj_{m,n}}{dm} dm > \int_m^{m+1} dm = 1. \tag{43}$$

We set $m = \operatorname{argmin}_{\mu \in \mathbb{N}_0} \{j_{\mu,1} \geq k_j R_0\}$, i.e., $m = M_-$. If we have that $j_{m,1} \in R_0(k_j - \Delta k, k_j + \Delta k) \subset (k_j R_0 - 1, k_j R_0 + 1)$, it follows from (43) that there can be no $j_{\tilde{m},1}$ with $\tilde{m} > m$ within the same interval. This implies that all the zeros $k_{\tilde{m},n} \in (k_j - \Delta k, k_j + \Delta k)$ must have order $\tilde{m} \leq m = M_-$ and hence that $\max_{m \in M_{k_j}} m \leq M_- \leq \beta_{k_j}$. $\square$

The result shows that, with a proper choice of $\Delta k$, we regularize the multi-frequency ISP by making a selective SVE [13] for each frequency, and it becomes clear that, under the additional requirement $\Delta k \leq 1/R_0$, the FB-coefficients of the source are stably recoverable with respect to noise in the measurements.

To recapitulate the developments of Theorem 1 and Theorem 2 in terms of reconstruction of sources, consider an ISP geometry specified by $R$ and $R_0$, as well as a finite-dimensional source space $S_{M,N}$. Assume $\Delta k$ is the largest number such that $\Delta k \leq 1/R_0$ and furthermore such that the inequality (27) holds for all indices $0 \leq m \leq M$ and $1 \leq i \leq N$.
If, for each $k_{m,i} \in Q_{M,N}$, there is some $k_j \in Q_r$ satisfying $|k_j - k_{m,i}| < \Delta k$, then any source $s \in S_{M,N}$ is uniquely recoverable from measurements at the frequencies $k_j \in Q_r$. Furthermore, the reconstruction is stable in the sense that it does not incorporate data associated with singular values with indices above the bandwidth $\beta_{k_j}$ for any of the measurements $U_{k_j}, k_j \in Q_r$.

We end this section by commenting that our analysis offers an explanation for a phenomenon observed by many authors working on the ISP; namely that the well-posedness of the ISP increases when the frequency increases. The results in Section 4.1 show that the density of the zeros of Bessel functions $J_m$ increases as one moves up the real line. Together with the results in Section 4.3, this implies that one can expect to stably approximate more SVE-coefficients from measurements at higher frequencies. Hence, in general, we are able to stably recover more information about the source from higher frequency measurements.

## 5 Numerical method and experiments

### 5.1 Numerical method

The numerical method consists of computing a reduced frequency set, implementing the $K$-matrix mapping FB-coefficients to SVE-coefficients, and computing the SVE-coefficients. We summarize the procedure below, and briefly elaborate on some of the steps.

1) Given a domain $D$, a source domain $D_0$, and wave speed $c$, one decides on a preferable resolution expressed in terms of the finite dimensional space $S_{M,N}$.

2) Compute the reduced frequency set $Q_s$ and assemble $Q_r$.

3) Use the frequencies from $Q_r$ to assemble the matrix $K$, according to equations (20)-(24).

4) Record the measurements $\{U_{k_j}\}_{k_j \in Q_s}$, where $k_j = \frac{\omega_j}{c}$ and $U_{k_j} = F_{k_j} s$. Use a quadrature method to compute the SVE-coefficients according to equation (42).

5) Solve $\hat{U} = K\hat{S}$ for $\hat{S}$, and express $s$ as $s = \sum_{i=-M}^{M} \sum_{j=1}^{N} \hat{s}_{i,j} \varphi_{i,j}$.

The problem specification in 1) dictates a choice of the measurement frequency set $Q_s$ in 2). From (34), we have available the formula for the entries of the matrix $K$, and we construct it using the entries in $Q_r$ according to (24)-(20).

In step 4) we sample each measurement $U_k$ on $\partial D$; the spatial sampling rate has to be small enough to avoid aliasing. Lastly, in 5) we express the solution as a finite linear combination of FB-vectors on a suitable mesh.

## 5.2 Numerical experiments

We test the method on three sources with different characteristics. The first source is a simple linear combination in $S_{3,3}$. Next, we take a smooth, slowly varying source, and last, a discontinuous source with some smaller, sharper features. In all cases we compute the reconstruction for different choices of $M$, $N$ and $\Delta k$, and also let $\Delta k$ violate the invertibility condition. The largest admissible $\Delta k$ from Lemma 1 is written in bold. The measurements are computed by solving equation with the built-in 2D-quadrature in MATLAB, and sampled at 200 uniformly distributed points on $\partial D$. We add complex Gaussian white noise to the measurements with relative noise level 0.2 (20%). The inner product in (42) is approximated using Simpson's rule. We compute a number of interesting $L^2$ error-norms for the reconstructions; they are listed in Tables 2-6. In the following, we denote the source by $s$, the reconstructed source by $s_r$, the source reconstructed from noisy measurements by $s_r^\varepsilon$, and the projection of $s$ on $S_{M,N}$ by $s_p$. In all reconstructions we let $D_0$ and $D$ have radii 1 and 1.5, respectively, and assume a wave speed $c = 1$.

### 5.2.1 Source in $S_{3,3}$

First, we consider the following linear combination of FB-basis functions:

$$s = 2\varphi_{0,1} + \pi \varphi_{3,3}. \tag{44}$$

We expect small errors, due only to the sampling and numerical approximations. This agrees with the observations in Table 2. The robustness towards noise, described in Section 4.3 is also evident, as the influence of the noise is seen to be very small.

### 5.2.2 Smooth source

Inspired by the papers [7, 26], we consider the source

$$s(x,y) = 0.3(1-3x)^2 \exp(-(3x)^2 - (3y+1)^2) - (0.2(3x) - (3x)^3$$
$$- (3y)^5) \exp(-(3x)^2 - (3y)^2) - 0.03 \exp(-(3x+1)^2 - (3y)^2).$$

We do reconstructions in $S_{3,3}$ and $S_{7,7}$, both with a number of reduced frequency sets corresponding to different choices of $\Delta k$. Results are presented in Tables 3 and 4, and plots of the source, the reconstruction, and the pointwise error are found in Figure 4. We observe that the we already get a quite accurate approximation of $s$ by just choosing $S_{3,3}$, and that we obtain fine reconstructions with as few as 5 frequencies. Also, all reconstructions are stable with respect to noise in the measurements. By comparing $\|s-s_r\|_{L^2}/\|s\|_{L^2}$ and $\|s-s_p\|_{L^2}/\|s\|_{L^2}$, we see that the effect of $s \notin S_{M,N}$ is small. Here, and later, we notice that $K$ is invertible, even for larger $\Delta k$, emphasizing that the invertibility condition in Theorem 1 is only a sufficient condition. For reconstructions in $S_{7,7}$, the projection becomes a more accurate approximation than the reconstruction, and we suspect this is due to the projection not suffering from the the approximation errors described in Section 3.3.

| $\|Q_s\|$ | 9 | 8 | 6 | **5** | 4 |
|---|---|---|---|---|---|
| $\Delta k$ | 0.25 | 0.5 | 0.75 | **0.91** | 1.5 |
| $\frac{\|s-s_r\|_{L^2}}{\|s\|_{L^2}}100\%$ | 1.3% | 1.3% | 1.4% | 1.40% | 1.3% |
| $\frac{\|s-s_r^\varepsilon\|_{L^2}}{\|s\|_{L^2}}100\%$ | 2.7% | 2.9% | 2.4% | 2.2% | 4.6% |

Table 2: Relative errors in percent, with and without noise, for a source $S_{3,3}$, measured in the $L^2$-norms.

| $\|Q_s\|$ | 9 | 8 | 6 | **5** | 4 |
|---|---|---|---|---|---|
| $\Delta k$ | 0.25 | 0.5 | 0.75 | **0.91** | 1.5 |
| $\frac{\|s-s_r\|_{L^2}}{\|s\|_{L^2}}100\%$ | 15.3% | 17.5% | 20.8% | 16.6% | 16.3% |
| $\frac{\|s-s_p\|_{L^2}}{\|s\|_{L^2}}100\%$ | 15.1% | 15.1% | 15.1% | 15.1% | 15.1% |
| $\frac{\|s_p-s_r\|_{L^2}}{\|s_p\|_{L^2}}100\%$ | 2.5% | 8.8% | 14.4% | 7.1% | 6.2% |
| $\frac{\|s-s_r^\varepsilon\|_{L^2}}{\|s\|_{L^2}}100\%$ | 15.5% | 18.1% | 21.5% | 17.4% | 16.4% |
| $\frac{\|s_p-s_r^\varepsilon\|_{L^2}}{\|s_p\|_{L^2}}100\%$ | 3.3% | 10.0% | 15.4% | 8.6% | 6.4% |

Table 3: Number of frequencies and relative errors in percent for the smooth source, measured in the $L^2$-norms, with $M = N = 3$.

### 5.2.3 Discontinuous source

Next, we consider the source

$$s_d(x,y) = 0.1 + \begin{cases} 1 & \text{for } (x,y) \in \hat{B}((-0.4,-0.08), 0.05), \\ 0.5 & \text{for } (x,y) \in \{(x,y): |x-0.2| \leq 0.15, |y+0.4| \leq 0.15\}, \\ 2 & \text{for } (x,y) \in \{(x,y): |x+0.2| \leq 0.2, |y-0.4| \leq 0.3\}. \end{cases} \quad (45)$$

We consider reconstructions in $S_{5,5}$ and $S_{15,15}$, also with different choices of reduced frequency set. Tables 5 and 6 contain error norms. A plot of the source, the reconstruction and the point-wise error is shown in Figure 5 and 6. It is evident from the plots that reconstructions in $S_{5,5}$ do not reproduce the smaller features, but that reconstructions in $S_{15,15}$ do quite well. We also notice the error around the edges, a common phenomenon when approximating with generalized Fourier-series. It is clear from the tables that the main error comes from the projection, as expected. Also, for reconstructions in $S_{15,15}$, the method shows significant robustness to noise in the measurements.

| $\lvert Q_s \rvert$ | 35 | **21** | 17 | 13 | 10 |
|---|---|---|---|---|---|
| $\Delta k$ | 0.25 | **0.61** | 0.75 | 1 | 1.5 |
| $\frac{\lVert s-s_r\rVert_{L^2}}{\lVert s\rVert_{L^2}}100\%$ | 6.5% | 7.8% | 19.2% | 19.0% | 20.7% |
| $\frac{\lVert s-s_p\rVert_{L^2}}{\lVert s\rVert_{L^2}}100\%$ | 1.6% | 1.6% | 1.6% | 1.6% | 1.6% |
| $\frac{\lVert s_p-s_r\rVert_{L^2}}{\lVert s_p\rVert_{L^2}}100\%$ | 6.3% | 7.6% | 19.2% | 18.9% | 20.6% |
| $\frac{\lVert s-s_r^\varepsilon\rVert_{L^2}}{\lVert s\rVert_{L^2}}100\%$ | 7.6% | 8.6% | 19.8% | 18.8% | 19.5% |
| $\frac{\lVert s_p-s_r^\varepsilon\rVert_{L^2}}{\lVert s_p\rVert_{L^2}}100\%$ | 7.4% | 8.5% | 19.7% | 18.7% | 19.4% |

Table 4: Number of frequencies and relative errors in percent for the smooth source, measured in the $L^2$-norms, with $M = N = 7$.

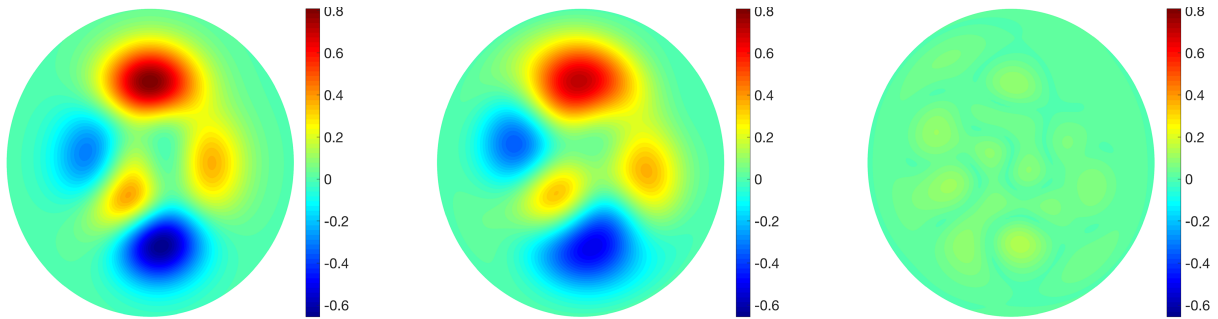

Figure 4: Smooth source. Left: Source. Middle: Reconstruction with $\Delta k = 0.91$, using 5 frequencies, and measurements with 20% additive Gaussian noise. Right: Pointwise absolute error.

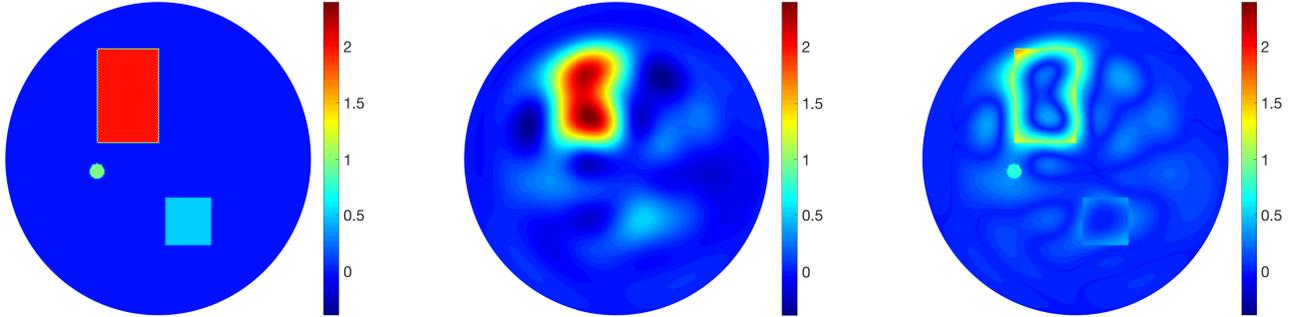

Figure 5: Discontinuous source. Left: Source. Middle: Reconstruction in $S_{5,5}$ with $\Delta k = 0.7$, using 14 frequencies, and measurements with 20% additive Gaussian noise. Right: Pointwise absolute error.

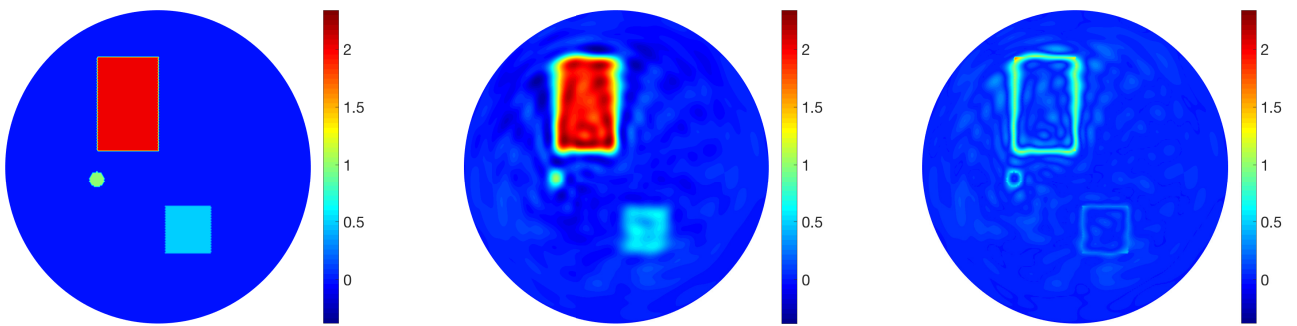

Figure 6: Discontinuous source. Left: Source. Middle: Reconstruction in $S_{15,15}$ with $\Delta k = 0.5$, using 56 frequencies, and data with 20% additive Gaussian noise. Right: Pointwise absolute error.

| $\lvert Q_s\rvert$ | 20 | **14** | 11 | 9 | 7 |
|---|---|---|---|---|---|
| $\Delta k$ | 0.25 | **0.7** | 0.75 | 1 | 1.5 |
| $\frac{\lVert s-s_r\rVert_{L^2}}{\lVert s\rVert_{L^2}}100\%$ | 42.4% | 43.2% | 42.9% | 43.4% | 43.8% |
| $\frac{\lVert s-s_p\rVert_{L^2}}{\lVert s\rVert_{L^2}}100\%$ | 42.2% | 42.2% | 42.2% | 42.2% | 42.2% |
| $\frac{\lVert s_p-s_r\rVert_{L^2}}{\lVert s_p\rVert_{L^2}}100\%$ | 3.2% | 9.9% | 8.6% | 11.3% | 12.6% |
| $\frac{\lVert s-s_r^\varepsilon\rVert_{L^2}}{\lVert s\rVert_{L^2}}100\%$ | 42.5% | 43.3% | 43.0% | 43.4% | 44.1% |
| $\frac{\lVert s_p-s_r^\varepsilon\rVert_{L^2}}{\lVert s_p\rVert_{L^2}}100\%$ | 4.2% | 10.2% | 9.1% | 11.2% | 14.1% |

Table 5: Relative errors in percent for the discontinuous source, measured in the $L^2$-norms, with $M = N = 5$.

| $\lvert Q_s\rvert$ | 91 | **56** | 40 | 30 | 22 |
|---|---|---|---|---|---|
| $\Delta k$ | 0.25 | **0.5** | 0.75 | 1 | 1.5 |
| $\frac{\lVert s-s_r\rVert_{L^2}}{\lVert s\rVert_{L^2}}100\%$ | 25.9% | 26.2% | 27.3% | 27.3% | 30.7% |
| $\frac{\lVert s-s_p\rVert_{L^2}}{\lVert s\rVert_{L^2}}100\%$ | 25.3% | 25.3% | 25.3% | 25.3% | 25.3% |
| $\frac{\lVert s_p-s_r\rVert_{L^2}}{\lVert s_p\rVert_{L^2}}100\%$ | 5.0% | 6.9% | 10.2% | 10.4% | 17.4% |
| $\frac{\lVert s-s_r^\varepsilon\rVert_{L^2}}{\lVert s\rVert_{L^2}}100\%$ | 26.3% | 26.6% | 27.5% | 27.9% | 30.1% |
| $\frac{\lVert s_p-s_r^\varepsilon\rVert_{L^2}}{\lVert s_p\rVert_{L^2}}100\%$ | 7.1% | 8.3% | 10.7% | 12.0% | 17.8% |

Table 6: Relative errors in percent for the discontinuous source, measured in the $L^2$-norms, with $M = N = 15$.

# 6 Conclusion

We have analyzed the multi-frequency inverse source problem with measurements at a finite number of frequencies. Our main results are the criteria which allow for a reduction of the frequency set while at the same time guaranteeing stable reconstructions of sources in finite-dimensional subspaces $S_{M,N} \subset L^2$. Further, we have devised a direct reconstruction scheme to reconstruct such sources. The method is stable with respect to noisy measurements, and regularization is implicit in the choice of the finite-dimensional source space and the corresponding measurement frequencies. Results from numerical experiments support our findings.

In the future, it would be interesting to see if our analysis can be generalized to more arbitrary geometries and higher dimensions, and to the Helmholtz equation with variable wave speed.

## Acknowledgements

Mirza Karamehmedović and Kim Knudsen would like to acknowledge support from The Danish Council for Independent Research — Natural Sciences.